\documentclass[11pt,letterpaper]{article}

\title{A Combinatorial Approach to the Groebner Bases for Ideals Generated by Elementary Symmetric Functions}
\author{AJ Bu}
\date{\today}

\usepackage[T1]{fontenc}
\usepackage{lmodern}
\usepackage[pdftex,margin=1in]{geometry}
\usepackage{enumerate}
\usepackage{amsmath}
\usepackage{amssymb}
\usepackage{fancyhdr}
\usepackage{lastpage}
\usepackage[alwaysadjust]{paralist}
\usepackage{amsthm}
\usepackage{graphicx}
\usepackage{wasysym}
\usepackage{xcolor}
\usepackage{listings}
\usepackage{hyperref}

\newcommand{\N}{\mathbb{N}}

\newcommand{\la}{\langle}
\newcommand{\ra}{\rangle}

\newtheorem{Thm}{Theorem}
\newtheorem{Prop}[Thm]{Proposition}

\newtheorem{Cor}[Thm]{Corollary}

\theoremstyle{definition}
\newtheorem{Def}[Thm]{Definition}

\theoremstyle{remark}

\theoremstyle{definition}

\theoremstyle{definition}

.

\newenvironment{Proof}{\noindent\textbf{Proof.}}{\qed}

\DeclareRobustCommand{\maketitle}{%
	\begingroup
	\begin{center}
	\Large \textbf{A Combinatorial Approach to the Groebner Bases for Ideals Generated by Elementary Symmetric Functions}\\[\medskipamount]
	\large AJ Bu \\[\medskipamount]
 \end{center}
	\endgroup
}

\begin{document}
	\maketitle
 \begin{abstract}
Previous work by Mora and Sala provides the reduced Groebner basis of the ideal  formed by the elementary symmetric polynomials in $n$ variables of degrees $k=1,\dots,n$, $\la e_{1,n}(x), \dots, e_{n,n}(x) \ra$.  Haglund, Rhoades, and Shimonozo expand upon this, finding the reduced Groebner basis of the ideal of elementary symmetric polynomials in $n$ variables of degree $d$ for $d=n-k+1,\dots,n$ for $k\leq n$.  In this paper, we further generalize their findings by using symbolic computation and experimentation to construct the reduced Groebner basis for the ideal generated by the elementary symmetric polynomials in $n$ variables of arbitrary degrees.
  \end{abstract}
  \noindent{\bf Keywords}: Groebner bases, symmetric functions
  
\noindent{\bf 2020 Mathematics Subject Classification}: 05E05

\section{Introduction}
In their paper \cite{ms}, Mora and Sala use computational and algebraic means to find the reduced Groebner basis of the ideal generated by the elementary symmetric polynomials in $n$ variables of degrees $d=1,\dots,n$. Haglund, Rhoades, and Shimonozo expand upon this, finding the reduced Groebner basis of the ideal of elementary symmetric polynomials in $n$ variables of degree $d$ for $d=n-k+1,\dots,n$ for $k\leq n$ \cite{HRS}.  In this paper, we further generalize their findings by using symbolic computation and experimentation to construct the Groebner basis for the ideal generated by the elementary symmetric polynomials in $n$ variables of arbitrary degrees.

\begin{Def} Let $k$ and $n$ be natural numbers. 
The \emph{elementary symmetric polynomial} of degree $k$ in $n$ variables, $x_1,\dots,x_n$, is
\[e_{k,n}(x) = \sum_{1 \leq i_1<\dots<i_k \leq n} x_{i_1} \dots x_{i_k}.\]
\end{Def}
\begin{Def}
The \emph{homogeneous symmetric polynomial} of degree $k$ in $n$ variables, $x_1,\dots,x_n$ is
\[h_{k,n}(x) = \sum_{1 \leq i_1\leq\dots\leq i_k \leq n} x_{i_1} \dots x_{i_k}.\]
\end{Def}
Given a set or multiset $S$ with elements in $\{1,\dots,n\},$ define the weight of $S$ to be 
\[wt(S) = \prod_{s \in S} x_s^{m(s)},\] where $m(s)$ is the multiplicity of $s$ in $S$.
For example,
\[wt(\{1,2,5\})=x_1 x_2 x_5, \text{ and } wt(\{1,1,3,4\})=x_1^2x_3x_4.\]
Then, $e_{k,n}(x)$ (respectively, $h_{k.n}(x)$) is the weight enumerator of the sets (respectively, multisets) with cardinality $k$ whose elements are in $\{1,\dots,n\}.$ Moreover, considering subsets of $\{1,\dots,n\}$ which do and do not contain $n$ separately, we have the following recursive definition
\begin{align*}    e_{k,n}(x) &= \begin{cases} 0, &\text{if } n<k\\
1, &\text{if } k=0\\
e_{k,n-1}(x)+x_n e_{k-1,n-1}(x), &\text{otherwise}.
\end{cases}
    \end{align*}
Similarly, when looking at multisets, we get
\begin{align*}    h_k(x_1,\dots,x_n) &= \begin{cases} 
0, &\text{if } n=0 \text{ and } k>0\\
1, &\text{if } k=0\\
h_{k,n-1}(x)+x_n h_{k-1,n}(x), &\text{otherwise}.
\end{cases} \end{align*}
We use these recursive definitions to write Maple functions $\texttt{eknS(x,k,n)}$ and $\texttt{hknS(x,k,n)}$, which output $e_{k.n}(x)$ and $h_{k,n}(x),$ respectively. These functions -- along with others used to investigate the Groebner basis of ideals generated by elementary symmetric polynomials -- can be found in the accompanying Maple package $\texttt{Solomon.txt},$ written by AJ Bu and Doron Zeilberger.

In \cite{ms}, Mora and Sala proved that $\{h_{1,n}(x),h_{2,n-1}(x),\dots,h_{n,1}(x) \}$ is a Groebner basis of the ideal $\la e_{1,n}(x),\dots,e_{n,n}(x) \ra$. Using the accompanying package to efficiently generate the reduced Groebner bases of many \emph{specific} ideals, we can extend their findings. We first use experimental methods to deduce a pattern for the reduced Groebner bases of the ideals $\la e_{1,n}(x),\dots,e_{k,n}(x)\ra $ and $\la e_{1,n}(x), e_{k,n}(x) \ra$ for arbitrary $k\leq n$, and prove them by combinatorial means.  We then investigate other cases to expand upon our results to the ideal $\la e_{k_1,n}(x), \dots, e_{k_m,n}(x) \ra$. We find a basis for this general case, proving that it generates the ideal, and show empirically that it is a Groebner basis.

\subsection{Groebner Bases}

\begin{Def}
A \emph{Groebner basis} of an ideal $I \subset k[x_1,...,x_n]$ (with respect to a given monomial order) is a finite subset $G=\{g_1,...,g_t\}$ of $I$ such that for that every nonzero polynomial $f$ in $I$, the leading term of $f$ is divisible by the leading term of $g_i$ for some $i$.

Moreover, it is  \emph{reduced} if, for every element $g \in G$, no monomial in $g$ is in $\la LT(G-{g})\ra,$ the ideal generated by the leading terms of the other elements in $G$.
\end{Def}

In order to efficiently determine whether or not a given basis is a Groebner basis, we use Buchberger's Criterion.

\begin{Thm}[Buchberger's Criterion]
	$G=\{g_1,\dots,g_t\}$ is a Groebner basis of $I$ with respect to a given monomial order if and only if $G$ generates $I$ and, for any distinct $g_i$ and $g_j$ in $G$, 
		\[\overline{S(g_i,g_j)}^G=0,\]
		where $\overline{S(f_i,f_j)}^G$ denotes the remainder of  the S-Polynomial of $f_i$ and $f_j$ upon division by $G$.
\end{Thm}
Given a polynomial $f$, we denote its leading monomial by $LM(f)$ and its leading term by $LT(f)$. The \emph{$S$-polynomial} of two polynomials $f$ and $g$ is
\[S(f,g) = \frac{LCM(LM(f),LM(g))}{LT(f)} f- \frac{LCM(LM(f),LM(g))}{LT(g)} g.\]

\begin{Thm}[The Division Algorithm in \text{$k[x_1,\dots,x_n]$}]
	Let $>$ be a fixed monomial order in $k[x_1,\dots,x_n]$. Let $F:=[f_1,\dots,f_m]$ be an ordered list of polynomials in $k[x_1,\dots,x_n]$. Then for any $f \in k[x_1,\dots,x_n]$, there exists $a_1,\dots,a_m, r \in k[x_1,\dots,x_n]$ such that
	 \\\indent 1.   $f= a_1f_1+\dots+a_m f_m +r$, 
	 \\\indent 2. for all $i$, either $a_if_i =0$ or $LT(f)\geq LT_{>} (a_if_i)$, and 
	 \\\indent 3. $r$ is a sum of monomials, none of which are divisible by any $LT(f_i)$.
	 \\We call r the \emph{remainder} of $f$ on $r$.
\end{Thm}
\label{sec:introduction}

\section{The Ideal \texorpdfstring{$\la e_{1,n}(x),\dots,e_{k,n}(x)\ra$}{TEXT}%
	} \label{sec:2}
The procedure $\texttt{Gkn(k,n,x)}$ in $\texttt{Solomon.txt}$ outputs the reduced Groebner basis (with respect to lexicographical order where $x_n > x_{n-1} > \dots > x_1$) for the ideal $\la e_{1,n}(x), \dots, e_{k,n}(x)\ra$. After running the procedure for multiple values $k$ and $n$, we can conjecture that the reduced Groebner basis is $\{h_{i,n-i+1}(x) | i=1\dots k\}$. Indeed, for the case $k=n$, this agrees with the Groebner basis that Mora and Sala proved in their paper \cite{ms}. In order to prove our conjecture, we use the following two relations between the elementary and homogeneous symmetric polynomials.

\subsection{Relations between elementary and homogeneous symmetric polynomials}
Recall that we defined the weight of $S=\{s_1,\dots,s_k\}$ with elements in $\{1,\dots,n\}$ to be
\[wt(S) := \prod_{s \in S} x_s^{m(s)},\]
where $m(s)$ is the multiplicity of $s$ in $S$.

\begin{Prop} \label{hkn} Let $k$ and $n$ be natural numbers. Then
	\[h_{k,n-k+1}(x) = \sum_{i=1}^{k} (-1)^{i+1} e_{i,n}(x) h_{k-i,n-k+1}(x)\]
\end{Prop}
\begin{Proof} This is equivalent to proving 
	\[\sum_{i=0}^{k} (-1)^{i} e_{i,n}(x) h_{k-i,n-k+1}(x)=0.\]
	This is trivial when $k>n$ because $h_{k-i,n-k+1}(x)=0$ when $0 \leq i \leq k-1$, and $e_{k,n}=0$. So, assume $k \leq n$.
	Then, the left-hand side is the weight enumerator of the set $\mathcal{S}_{k,n}$ of pairs $(A,B)$, where
	\begin{itemize}
	    \item $A$ is a subset of $\{1,\dots,n\}$ of order $|A|$,
	    \item $B$ is a multiset with cardinality $k-|A|$ whose elements are in $\{1,\dots,n-k+1\}$,
	\end{itemize}
	and the weight of $(A,B)$ is 
	\[w(A,B) = (-1)^{|A|} wt(A) wt(B).\]
	
	Let $f: \mathcal{S}_{k,n} \to \mathcal{S}_{k,n}$ be defined as 
	\[f(A,B) = \begin{cases} (A \cup \{\min(B)\}, B-\{\min(B)\}) , & \text{if } \min(B)<\min(A) \\
	(A \smallsetminus \{\min(A)\}, B +\{\min(A)\}), &\text{otherwise}.
	\end{cases}\]
	Note that this mapping is defined for all possible pairs of sequences, and it changes sign since the size of the first subset is either increasing or decreasing by $1$. 
	Moreover, if
	$\min(A) > \min(B)$ then 
	\begin{align*}
	    f(A,B) &=  (A \cup \{\min(B)\}, B-\{\min(B)\}) =: (A',B'), \text{ and}\\
	    f(A',B') &=(A,B),
	\end{align*}
	since clearly $\min(A')=\min(B) \leq \min(B')$.
	If $\min(A) \leq \min(B)$ then 
	\begin{align*}
	    f(A,B) &= (A \smallsetminus \{\min(A)\}, B +\{\min(A)\} ) =: (A',B'), \text{ and}\\
	    f(A',B') &=(A,B),
	\end{align*} 
	since $\min(B')=\min(A) < \min(A')$. Thus, all elements of $\mathcal{S}_{k,n}$ can be paired up into mutually cancelling pairs, concluding our proof. \end{Proof}

\begin{Prop} \label{ekn} For any $n,k \in \N$, 
	\[e_{k,n}(x) = \sum_{i=1}^{k} (-1)^{i+1} h_{i,n-i+1}(x)e_{k-i,n-i}(x).\]
\end{Prop}
\begin{Proof}
	Note that this is equivalent to 
	\[ \sum_{i=0}^{k} (-1)^{i} h_{i,n-i+1}(x)e_{k-i,n-i}(x)=0.\]
	
Again, this is trivial for $k>n$, so assume $k \leq n$.	The left-hand side is the weight enumerator of the set $\mathcal{S}_{k,n}$ of ordered pairs $(A,B)$, where
	\begin{itemize}
	    \item $A$ is a multiset with elements in $\{1,\dots, n-|A|+1\},$ where $|A|$ is the cardinality of $A$,
	    \item $B$ is a subset of $\{1,\dots, n-|A| \}$ of order $|B|:=k-|A|,$
	\end{itemize}
	and the weight of $(A,B)$ is 
	\[w(A,B) = (-1)^{|A|} wt(A) wt(B).\]
		Let $f: \mathcal{S}_{k,n} \to \mathcal{S}_{k,n}$ be defined as 
	\[f(A,B) = \begin{cases} (A + \{\min(B)\}, B \smallsetminus \{\min(B)\}) , & \text{if } \min(B)<\min(A) \\
	(A - \{\min(A)\}, B \cup \{\min(A)\}), &\text{otherwise}.
	\end{cases}\]
	As in the previous proof,  this involution pairs all elements of $\mathcal{S}_{k,n}$ into mutually cancelling pairs.
	\end{Proof}

\subsection{The Reduced Groebner Basis of the Ideal}
\begin{Prop} \label{ekn, hkn ideal} Let $k$ and $n$ be natural numbers such that $k \leq n$.
\[\la e_{1,n}(x),\dots,e_{k,n}(x)\ra =\la h_{1,n}(x), h_{2,n-1}(x),\dots , h_{k,n-k+1}(x)\ra.\]
\end{Prop}
\begin{Proof}
	For $i=1,\dots,k$, we have 
	\begin{align*}
	h_{i,n-i+1}(x) &\in \la e_{1,n}(x),\dots,e_{k,n}(x) \ra, &\text{ and }&&	e_{i,n}(x) &\in \la h_{1,n}(x), h_{2,n-1}(x),\dots, h_{k,n-k+1}(x)\ra 
	\end{align*}  by propositions \ref{hkn} and \ref{ekn}, respectively.  It immediately follows that \[\la e_{1,n}(x),\dots,e_{k,n}(x)\ra =\la h_{1,n}(x), h_{2,n-1}(x),\dots , h_{k,n-k+1}(x)\ra.\]
\end{Proof}

\begin{Prop} \label{ekn, hkn ideal groeb} Let $k$ and $n$ be natural numbers.  The set $G:=\{ h_{i,n-i+1}(x) \mid 1 \leq i \leq k \}$ is the reduced Groebner basis of the ideal $\la e_{1,n}(x),\dots,e_{k,n}(x)\ra$ with respect to lexicographical order, where $x_n>x_{n-1}>\dots>x_1$.
	\end{Prop}
\begin{Proof}
	By proposition \ref{ekn, hkn ideal}, the set $G$ generates the ideal $I:=\la e_{1,n}(x),\dots,e_{k,n}(x)\ra $. The $S$-polynomial of any two distinct elements $h_{i,n-i+1}(x)$ and $h_{j,n-j+1}(x)$ in $G$ is
	\begin{align*}
	S(h_{i,n-i+1}(x),h_{j,n-j+1}(x)) &=x_{n-j+1}^j h_{i,n-i+1}(x)-x_{n-i+1}^i h_{j,n-j+1}(x)
	\\&=h_{j,n-j+1}(x)\sum_{\ell=0}^{i-1} x_{n-i+1}^\ell h_{i-\ell,n-i}(x)-h_{i, n-i+1}(x)\sum_{\ell=0}^{j-1} x_{n-j+1}^\ell h_{j-\ell,n-j}(x).
	\end{align*}
	To prove the second equality, note that it is equivalent to 
		\[h_{i, n-i+1}(x)\sum_{\ell=0}^{j} x_{n-j+1}^\ell h_{j-\ell,n-j} (x)=   h_{j,n-j+1}(x)\sum_{\ell=0}^{i} x_{n-i+1}^\ell h_{i-\ell, n-i}(x).\]
	$x_{n-j+1}^\ell h_{j-\ell,n-j} (x)$ is the weight enumerator of all multisets of cardinality $j$ with elements taken from $\{1,\dots, n-j+1\}$, where $n-j+1$ appears exactly $\ell$ times. Thus, it is clear that 
	\[\sum_{\ell=0}^{j} x_{n-j+1}^\ell h_{j-\ell,n-j} (x) =h_{j,n-j+1}(x).\]
	It follows that
		\begin{align*}
	S(h_{i,n-i+1}(x),h_{j,n-j+1}(x)) &=h_{j,n-j+1}(x)\sum_{\ell=0}^{i-1} x_{n-i+1}^\ell h_{i-\ell,n-i}(x)-h_{i, n-i+1}(x)\sum_{\ell=0}^{j-1} x_{n-j+1}^\ell h_{j-\ell,n-j}(x).
	\end{align*}
	Moreover, for $i \neq j$
	\begin{align*}
	LT\left(h_{i, n-j+1}(x) \sum_{\ell=0}^{i-1} x_{n-i+1}^\ell h_{i-\ell,n-i}(x) \right)  &= x_{n-j+1}^j x_{n-i+1}^{i-1} x_{n-i}\\
	& \neq   x_{n-i+1}^i x_{n-j+1}^{j-1} x_{n-j}\\
	&= LT \left(h_{i,n-i+1}(x)\sum_{\ell=0}^{j-1} x_{n-j+1}^\ell h_{j-\ell,n-j}(x)\right).
	\end{align*}
	Hence, 
	\[LT(	S(h_{i,n-i+1}(x),h_{j,n-j+1}(x)))= \max( x_{n-j+1}^j x_{n-i+1}^{i-1} x_{n-i},  x_{n-i+1}^i x_{n-j+1}^{j-1} x_{n-j} )\]
	and, by the division algorithm,  
	\[	\overline{S(h_{i,n-i+1}(x),h_{j,n-j+1}(x))}^G=0.\]
	Therefore, $G$ is a Groebner basis of $I$ by Buchberger's Criterion. Furthermore, $G$ is a reduced Groebner basis because, for any distinct $i,j$, $LT(h_{i,n-i+1}(x))=x_{n-i+1}^i$ cannot divide the terms in  $h_{j,n-j+1}(x)$. This follows from the fact that the elements of $h_{j,n-j+1}(x)$ have lower degree if $i>j$, and they cannot be multiples of $x_{n-i+1}$ if $i<j$.
\end{Proof}

\section{The Ideal \texorpdfstring{$\la e_{1,n}(x),e_{k,n}(x)\ra$}{TEXT}%
	}

\begin{Prop} \label{e1 ek} Let $k$ and $n$ be natural numbers such that $k \leq n$. Then
	\[\la e_{1,n}(x), e_{k,n}(x)\ra = \la e_{1,n}(x),e_{1,n-1}(x) e_{k-1,n-1}(x)-e_{k,n-1}(x)\ra. \]
\end{Prop}
\begin{Proof} It suffices to prove
	\begin{align*}
		e_{1,n-1}(x) e_{k-1,n-1}(x)-e_{k,n-1}(x) &= e_{1,n}(x)e_{k-1,n-1}(x)-e_{k,n}(x),
	\end{align*}
	or equivalently
	\[e_{k,n}(x)-e_{k,n-1}(x) = x_ne_{k-1,n-1}(x).\]
	This clearly holds since both sides are the weight enumerator of the $k$ element subsets of $\{x_1,\dots,x_n\}$ containing $x_n$.
\end{Proof}

\begin{Prop} \label{e1 ek groeb} Let $k$ and $n$ be natural numbers such that $n \geq k$. Then
\[G:=\{ e_{1,n}(x),e_{1,n-1}(x) e_{k-1,n-1}(x)-e_{k,n-1}(x) \}\] is the reduced Groebner basis of the ideal $\la e_{1,n}(x), e_{k,n}(x) \ra$ with respect to lexicographical order, where $x_n>x_{n-1}>\dots>x_1$.
\end{Prop}
\begin{Proof}
	By proposition \ref{e1 ek}, $G$ generates $I:=\la e_{1,n}(x), e_{k,n}(x) \ra $. Taking the $S$-polynomial of the elements in $G$, we have
	\begin{align*}
	S(e_{1,n}(x),&e_{1,n-1}(x) e_{k-1,n-1}(x)-e_{k,n-1}(x))
	\\&=x_{n-1}^2x_{n-2}\dots x_{n-k+1}e_{1,n}(x)-x_n \big( e_{1,n-1}(x) e_{k-1,n-1}(x)-e_{k,n-1}(x) \big)\\
	&=\big(e_{1,n-1}(x) e_{k-1,n-1}(x)-e_{k,n-1}(x)\big) e_{1,n-1}(x)
	\\&\hspace*{25pt} -e_{1,n}(x)\big(e_{1,n-1}(x) e_{k-1,n-1}(x)-e_{k,n-1}(x)- x_{n-1}^2x_{n-2}\dots x_{n-k+1}\big)
	\end{align*}
	Note that the second equality obviously holds since it can be rewritten as
	\[e_{1,n}(x)\big(e_{1,n-1}(x) e_{k-1,n-1}(x)-e_{k,n-1}(x)) =\big(e_{1,n-1}(x) e_{k-1,n-1}(x)-e_{k,n-1}(x))e_{1,n}(x).\]
	It is also clear that 
\begin{align*}
LT \big( (e_{1,n-1}(x)& e_{k-1,n-1}(x)-e_{k,n-1}(x)) e_{1,n-1}(x)\big) \\&< LT\big( e_{1,n}(x)(e_{1,n-1}(x) e_{k-1,n-1}(x)-e_{k,n-1}(x)- x_{n-1}^2x_{n-2}\dots x_{n-k+1})\big) ,
\end{align*}
since the latter is a multiple of $x_n$. Hence, 
\[\overline{S(e_{1,n}(x),e_{1,n-1}(x) e_{k-1,n-1}(x)-e_{k,n-1}(x))}^G=0, \] and $G$ is a Groebner basis. It is clearly reduced since no term in $e_{1,n-1}(x) e_{k-1,n-1}(x)-e_{k,n-1}(x)$ is divisible by $x_n$ and no term in $e_{1,n}(x)$ is divisible by $x_{n-1}^2$.
\end{Proof}

\section{Investigation into the General Case}
The procedure $\texttt{GSn(S,n,x)}$ in $\texttt{Solomon.txt}$ inputs a set $S=\{k_1,\dots,k_m\}$, a non-negative integer $n$, and a variable $x$. It outputs the reduced Groebner basis (with respect to lexicographical order where $x_n > x_{n-1} > \dots > x_1$ ) for the ideal $\la e_{k_1,n}(x), \dots, e_{k_m,n}(x)\ra$. Using this procedure to analyze the reduced Groebner bases for various ideals, we conjecture the following basis for arbitrary $S$ and $n$.

\begin{Prop} \label{gen} Let $k_1,\dots,k_m$, and $n$ be positive integers such that $1 \leq k_1<\dots<k_m \leq n.$ Let $I$ be the ideal 
$$I:=\la e_{k_1,n}(x), \dots , e_{k_m,n}(x)\ra,$$
and let $M$ be the set of matrices of the form
$$\begin{bmatrix}
	e_{k_m- i_{m-1}, n-i_{m-1} } (x) & \dots & e_{k_m- i_{1}, n-i_{1} } (x)  & e_{k_m , n} (x)  \\
	\vdots & \dots & &\vdots\\
e_{k_1- i_{m-1}, n-i_{m-1} } (x) & \dots & e_{k_1- i_{1}, n-i_{1} } (x)  & e_{k_1 , n} (x)
	\end{bmatrix},$$
where $i_1 \in \{ 1, 2, \dots,  k_1, k_2  \}$ and 
	$	i_j \in \{ i_{j-1}+1,i_{j-1}+2, \dots, k_j, k_{j+1}\}$ for  $j>1$.
Then the set 
$$G:=\{ det(m) \mid m \in M\}$$
is a basis of $I$.
\end{Prop}

\begin{Proof} 
Note that the entries of the last column of any matrix in $M$ are the elementary symmetric polynomials
$$e_{k_1,n}(x), \dots, e_{k_m,n}(x)$$
 that generate $I$. It immediately follows that $\la G \ra \subseteq I$.

For the other containment, let $m_1$ be the matrix in $M$ where $i_j=k_{j+1}$. Then,
\begin{align*}
     \det(m_1)&=\det  \left( 
\begin{bmatrix}
1 & e_{k_m-k_{m-1},n-k_{m-1}}(x) & \dots & e_{k_m-k_2,n-k_2}(x) & e_{k_m,n}(x)\\
0 & 1 & \dots &  e_{k_{m-1}-k_2,n-k_2}(x)& e_{k_{m-1},n}(x)\\
\vdots & 0 &\ddots &  \vdots& \vdots\\
\vdots & \vdots & \ddots & 1 &  e_{k_2,n}(x)\\
0 & 0 & \dots & 0 & e_{k_1,n}(x)\\
\end{bmatrix}
\right)\\
&=e_{k_1,n}(x).
\end{align*}
Therefore, $e_{k_1,n}(x) \in \la G \ra$.  Now suppose that for $L>1$, $e_{k_\ell,n}(x) \in \la G \ra$ for all $1\leq \ell<L.$ Let $m_L$ denote the matrix in $M$ such that  $i_j=k_j$ for $j<L$ and $i_j= k_{j+1}$ for $j \geq L$. Then, 
$$
m_{L} = \begin{bmatrix}
A_{L} & B_{L} \\
0 & C_{L}
\end{bmatrix},$$
where $A$ is an $(m-L) \times (m-L)$ triangular matrix with whose diagonal entries are all $1$, and $0$ is an  $ L \times (m-L)$ zero matrix. Therefore,
$\det m_L= \det C_L,$ where $C_L$ is the $L \times L$ matrix
$$
\begin{bmatrix}
 e_{k_{L} -k_{L-1},n-k_{L-1}}(x) &e_{k_{L}-k_{L-2},n-k_{L-2}}(x) &\dots& e_{k_L-k_1,n-k_1}(x)& e_{k_{L},n}(x)\\
1&   e_{k_{L-1}-k_{L-2},n-k_{L-2}}(x) &\dots & e_{k_{L-1}-k_1,n-k_1}(x) & e_{k_{L-1},n}(x)\\
0&1   & \ddots & \vdots& \vdots\\
\vdots &\ddots&\ddots&e_{k_{2}-k_1,n-k_1}(x) &  e_{k_2,n}(x)\\
 0&\dots  &0& 1 & e_{k_1,n}(x)
\end{bmatrix}.
$$
Define $c_{i}$ to be the $(L-1) \times (L-1)$ matrix formed by removing the $i-th$ row and the last column from $C_L$. Then,
\begin{align*}
\det(C_L) &= \sum_{i=1}^{L} (-1)^{i+1} e_{k_{i},n}(x) \det(c_{L+1-i}).
\end{align*}
Since $\det c_1=1$, it follows that
\begin{align*}
    \det M &= \det C_L\\
    &= (-1)^{L+1} e_{k_L,n}(x)+\sum_{i=1}^{L-1} (-1)^{i+1} e_{k_{i},n}(x) \det(c_{L+1-i}).
\end{align*}
Since $\det M \in \la G \ra $ and, by our inductive hypothesis, $e_{k_i,n}(x) \in \la G \ra$ for $1 \leq i \leq L-1,$ 
it follows $e_{k_L,n}(x)$ is in the ideal as well.  Thus, $e_{k_i,n}(x)\in \la G \ra $ for $i=1,\dots, m$, and  $I = \la G\ra$.
\end{Proof}
\vspace{11pt}

Since we found this basis by studying specifically the reduced Groebner bases of various ideals, we further conjecture that it is the reduced Groebner basis for arbitrary $S$ and $n$.  Indeed, it is obvious that the reduced Groebner basis for the ideal $\la e_{1,n}(x), e_{k,n}(x) \ra$ given in proposition \ref{e1 ek groeb} is of this form. Furthermore, upon further investigation, one can see that this Groebner basis for the general case also gives the reduced Groebner basis of the ideal $\la e_{1,n}(x),\dots,e_{k,n}(x)\ra$, which is stated in proposition \ref{ekn, hkn ideal groeb}. We show that this is true in the following proposition.

\begin{Prop}
 Let $k$ and $n$ be positive integers such that $1 \leq k \leq n.$ Let $I$ be the ideal 
$$I:=\la e_{1,n}(x), \dots , e_{k,n}(x)\ra,$$
and let $M$ be the set of matrices of the form
$$\begin{bmatrix}
	e_{k- i_{k-1}, n-i_{k-1} } (x) & \dots & e_{k- i_{1}, n-i_{1} } (x)  & e_{k , n} (x)  \\
	\vdots & \dots & \vdots &\vdots\\
e_{1- i_{k-1}, n-i_{k-1} } (x) & \dots & e_{1- i_{1}, n-i_{1} } (x)  & e_{1 , n} (x)
	\end{bmatrix},$$
where $1\leq i_1 < \dots < i_{k-1} \leq k$.
Then the set 
$$G:=\{ det(m) \mid m \in M\}$$
is the reduced Groebner basis of $I$.
\end{Prop}
\begin{Proof}
By the previous proposition, $G$ is a basis of $I$. Moreover, by proposition \ref{ekn, hkn ideal groeb}, the reduced Groebner basis of $I$ is $G':=\{h_{i,n-i}(x) \mid 1 \leq i \leq k \}$. Thus, it suffices to prove that $G=G'$.  

For any positive integer $L$, let $m_L$ denote the matrix such that no $i_j=L$.  Then, as shown in the previous proof, 
\begin{align*}
    \det m_L &=\det C_{L,n},
\end{align*}
where $$
C_{L,n}=\begin{bmatrix}
 e_{1,n-L+1}(x) &e_{2,n-L+2(x)} &\dots& e_{L-1,n-1}(x)& e_{L,n}(x)\\
1&   e_{1,n-L+2}(x) &\dots & e_{L-2,n-1}(x) & e_{L-1,n}(x)\\
0&1   & \ddots & \vdots& \vdots\\
\vdots &\ddots&\ddots&e_{1,n-1}(x) &  e_{2,n}(x)\\
 0&\dots  &0& 1 & e_{1,n}(x)
\end{bmatrix}.$$

For any positive integers $L$ and $n$ where $L \leq n$, we claim $$\det C_{L,n}=h_{L,n-L+1}(x).$$ To prove this claim, we use induction over $L$, where for each $L\in\mathbb{N}$, we will show that the claim holds for all $n$.  We begin with the base case; for any positive integer $n$, we have $C_{1,n}= \begin{bmatrix}
 e_{1,n}(x)
\end{bmatrix},$ so clearly 
\begin{align*}
    \det C_{1,n} &=
 e_{1,n}(x)
\\
&=h_{1,n}(x).
\end{align*}

Now suppose that for any given $L\geq 2$, we have $\det C_{\ell,N }=h_{\ell,N-\ell+1}(x)$ for any $0<\ell<L$ and $N>\ell.$ Since
$$\det C_{L,n} =\sum_{i=1}^{L} (-1)^{i+1} e_{i,n}(x) \det(c_{L+1-i}) ,$$
where $c_i$ is formed by removing the $i-th$ row and the last column from $C_{L,n},$ and we have shown that 
	\[h_{k,n-k+1}(x) = \sum_{i=1}^{k} (-1)^{i+1} e_{i,n}(x) h_{k-i,n-k+1}(x), \]
	it is enough to show that $\det c_i = h_{i-1,n-L+1}.$ Since $c_1$ is a triangular matrix whose diagonal entries are $1$, it is obvious that 
	\begin{align*} \det c_1 & = 1
	\\ &= h_{0, n-L+1}.\end{align*}
	
	For $i>1$, $c_i$ can be written as
	$$c_i=\begin{bmatrix}
	 a_i & b_i \\0& d_i,
	\end{bmatrix}$$
	where $d_i$ is an $(L-i) \times (L-i)$ triangular matrix whose diagonal entries are all $1$,  and $a_i=C_{i-1,n-L+i-1}$. Therefore,
	\begin{align*}
	    \det c_i &= \det C_{i-1,n-L+i-1}\\
	    &= h_{i-1,n-L+1},
	\end{align*}
    by our inductive hypothesis. Thus,
    \begin{align*}
    \det C_{L,n} &= \sum_{i=1}^{L} (-1)^{i+1} e_{i,n}(x) \det(c_{L+1-i})\\
     &=\sum_{i=1}^{L} (-1)^{i+1} e_{i,n}(x) h_{L-i,n-L+1}(x)\\
     &= h_{L,n-L+1}(x),
\end{align*}
as desired.
\end{Proof}

\section{Conclusion}
In section \ref{sec:2}, we showed that the reduced Groebner basis of the ideal $\la e_{1,k}(x),\dots, e_{k,n}(x) \ra$ is $\{h_{i,n-i+1}(x)| i=1, \dots, k \}$. The following corollary to this proposition was proved by Chevally \cite{ch}.
\begin{Cor}
The Hilbert series of $K[x_1,\dots,x_n]/\la e_1(x), \dots, e_k (x) \ra$ is 
\[\frac{\prod_{i=1}^n (1-t^i)}{(1-t)^n}.\]
In particular, the dimension of this finite dimensional vector space is $n!$.
\end{Cor}
    
We have given a formal proof for a basis of $\la e_{k_1,n}(x), \dots,e_{k_m,n}(x) \ra $,  the ideal generated by an arbitrary set of elementary symmetric functions of degree $n$, and we proved that it is the reduced Groebner basis for ideals of the form $\la e_{1,n}(x), e_{k,n}(x) \ra$ and $\la e_{1,n}(x), \dots, e_{k,n}(x) \ra.$  Using the procedure $\texttt{CheckConjGSn(S,n,x)}$ in $\texttt{Solomon.txt},$ we can verify that the generalized basis given in proposition \ref{gen} is the reduced Groebner basis of the ideal $\la e_{k_1,n}(x),\dots, e_{k_m}(x)\ra$ for a given set $S=\{k_1,\dots,k_m\}$ and positive integer $n$.  Verifying that this is the case for many $S$ and $n$, we can empirically show that this basis is the reduced Groebner basis for the ideal generated by an arbitrary set of elementary symmetric functions.

One direction for further research is to formally prove that the basis we have found for the general case is the reduced Groebner basis.  We can also try to find similar identities for other ideals, such as those generated by various power sum symmetric polynomials or homogeneous symmetric polynomials of arbitrary degrees.

\vspace{11pt}
\noindent \textbf{Acknowledgement:} Thank you to Doron Zeilberger for helpful feedback and guidance in research for this paper. Also, thank you to Richard Stanley for helpful feedback regarding past results relevant to this paper. Thank you to Brendon Rhoades for reaching out regarding the relevant findings in $\cite{HRS}.$


\begin{thebibliography}{1}
\bibitem{B} F. Bergeron. Algebraic Combinatorics and Coinvariant Spaces. CMS Treatises in Mathematics. Boca Raton:
Taylor and Francis, 2009

\bibitem{ch} C. Chevalley, \textit{Invariates of finite groups generated by reflections}, Amer. J. Math, \textbf{67} (1955), 778-782.

\bibitem{CLO} D. Cox, J. B. Little, and D. O'Shea, "Using Algebraic Geometry," Springer, 1998.

\bibitem{HRS} J. Haglund, B. Rhoades, and M. Shimozono, \textit{Ordered set partitions, generalized coinvariant algebras, and the delta conjecture}, preprint, arXiv:\href{https://arxiv.org/pdf/1609.07575.pdf}{1609.07575v4} (2019).

\bibitem{ms} T. Mora and M. Sala, \textit{On the Gr\"{o}bner bases of some symmetric systems and their application to coding theory}, J. Symbolic Comput. \textbf{35} (2003) no. 2, 177-194.

\bibitem{S} B. Sturmfels. Algorithms in Invariant Theory. Springer-Verlag, Berlin, 1993.

\bibitem{z} D. Zeilberger, \emph{A combinatorial proof of Newton's identities}, Discrete Math. \textbf{49} (1984), 319.
\end{thebibliography}
\end{document}